\documentclass[12pt,twoside]{article}
\usepackage{amsfonts}
\usepackage{amssymb}
\usepackage{amsmath}
{\markboth{\sl S.
Hamem and L. Kamoun}{\sl Uncertainty Principle Inequalities Related to Laguerre-Bessel Transform} \thispagestyle{plain}
\newtheorem{prop}{Proposition}[section]
\newtheorem{thm}{Theorem}[section]

\newtheorem{lem}{Lemma}[section]

\begin{document}
\title{Uncertainty Principle Inequalities Related to Laguerre-Bessel Transform\footnote{\small The authors are supported by the DGRST research project 04/UR/15-02}}
\author{Soumeya Hamem$^1$  and Lotfi Kamoun$^2$ \\\\\small Department of Mathematics,
Faculty of Sciences of Monastir,\\ \small University of Monastir, 5019 Monastir, Tunisia\\ \small $^1$E-mail : soumeya.hamem@yahoo.fr \\  \small
$^2$ E-mail : kamoun.lotfi@planet.tn }
\date{}
\maketitle
\begin{abstract} \noindent
In this paper, an analogous of Heisenberg inequality is established for Laguerre-Bessel transform. Also, a local uncertainty principle for this transform is investigated.
\end{abstract}
\section{Introduction}
The uncertainty principle states that a nonzero function and its Fourier transform cannot both be sharply localized.
In the language of quantum mechanics, this principle says that an observer cannot simultaneously and precisely determines the values of position and
momentum of a quantum particule. A mathematical formulation of this physical ideas is firstly developed by Heisenberg \cite{Hei} in $1927$.
For $f \in L^2(\mathbb{R}),$ a precise quantitative formulation of the uncertainty principle, usually called Heisenberg inequality, is the following
\begin{equation}\label{Hei}
\int_{\mathbb{R}}x^2\left|f(x)\right|^2dx.\int_{\mathbb{R}}\xi^2\left|\widehat{f}(\xi)\right|^2d\xi\geq\frac{1}{4}\left(\displaystyle\int_{\mathbb{R}}
\left|f(x)\right|^2dx\right)^2
\end{equation}
where $$\widehat{f}(\xi)=\displaystyle\frac{1}{\sqrt{2\pi}}\ \displaystyle\int_{\mathbb{R}}f(x)e^{-i\xi x}dx.$$
This result does not appear in Heisenberg paper \cite{Hei}. The relation (\ref{Hei}) appears in Weyl \cite{Wey} who credits the result to Pauli.
In framework of Hankel transform, Bowie in \cite{Bow} studied the Heisenberg uncertainty principle. R$\ddot{o}$sler in \cite{Ros} and Shimeno in \cite{Shi} have proved, by different methods, an Heisenberg inequality for the Dunkl transform. Recently, Ma in \cite{Rui} has obtained an Heisenberg inequality for the Jacobi transform.
Since the $20$'s of last century, many works have been devoted to studyng uncertainty principle in various forms. Among these, we can cite the works of  Faris \cite{far} and Price (\cite{pri}, \cite{pric}), whose aim is to establish local uncertainty inequalities. In this paper, firstly we obtain an analogous of Heisenberg inequality for the Laguerre-Bessel transform. Next, for this transform we develop further inequalities in the sharpest forms, which constitue the principle of local uncertainty. Throughout the paper, we denote $\mathbb{K}=[0,+\infty)\times [0,+\infty)$, $\widehat{\mathbb{K}}= [0,+\infty)\times \mathbb{N}$ and we designate by $C$ a positive constant, which is not necessarily same at each occurrence.
\section{Laguerre-Bessel transform}
In this section, we collect some notations and results about the Laguerre-Bessel harmonics analysis.
For more details, we refer the reader to \cite{Solta}. \\
For $\alpha\geq0$, we consider the following system of partial differential operators
$$\left\{
\begin{array}{lll}D_1&=&\displaystyle\frac{\partial^2}{\partial t^2}+\frac{2\alpha}{t}\frac{\partial}{\partial t}\\ &&\\
 D_2&=&\displaystyle\frac{\partial^2}{\partial x^2}+\frac{2\alpha+1}{x}\frac{\partial}{\partial x}+x^2D_1\,,\quad (x,t)\in \mathbb{K}.\end{array}\right.$$
For $(\lambda,m)\in \widehat{\mathbb{K}}$, the system
$$\left\{\begin{array}{l} D_1u=-\lambda^2 u\\ \\ D_2u=-2\lambda(2m+\alpha+1)u\\ \\ u(0,0)=1\,,\ \ \displaystyle\frac{\partial u}{\partial x}(0,0)
=\frac{\partial u}{\partial t}(0,0)=0\end{array}\right.$$
possesses a unique solution denoted $\varphi_{(\lambda,m)}$ and given by $$\varphi_{(\lambda,m)}(x,t)=j_{\alpha-\frac{1}{2}}(\lambda t)\
\mathcal{L}_m^\alpha(\lambda x^2),\ \ \ \ \ \ \  (x,t)\in \mathbb{K},$$
where $j_\alpha$ is the normalized Bessel function given by $$j_\alpha(x)=\Gamma(\alpha+1)\displaystyle\sum_{k=0}^{\infty}
\displaystyle\frac{(-1)^k}{k!\ \Gamma(\alpha+k+1)}\left(\displaystyle\frac{x}{2}\right)^{2k}$$
and $\mathcal{L}_m^\alpha$ is the Laguerre function defined on $[0,+\infty)$ by $$\mathcal{L}_m^{\alpha}(x)=
\displaystyle \frac{e^{-\frac{x}{2}}L_m^{\alpha}(x)}{L_m^{\alpha}(0)}$$
$L_m^{\alpha}$ being the Laguerre polynomial of degree $m$
and order $\alpha$ given by $$ L_m^\alpha(x)=\displaystyle\sum_{j=0}^{m}\frac{\Gamma(m+\alpha+1)(-x)^j}{\Gamma(m-j+1)\Gamma(j+\alpha+1)j!}.$$
$L_m^\alpha$ defined in terms of the generating function by
\begin{equation}\label{lag}
 \displaystyle\sum_{m=0}^{+\infty}t^m \ L_m^\alpha (x)=\displaystyle\frac{1}{(1-t)^{\alpha+1}}\ e^{-\frac{xt}{1-t}}
\end{equation}
 $\textbf{Notations}$ \\\\
 $\bullet$ $\mathcal{S}_*(\mathbb{K})$ the space of $C^\infty$ functions on $\mathbb{R}^2$, even with respect to each variable and rapidly decreasing
  together with all their derivatives i.e for all $k,p,q\ \in\mathbb{N}$,$$N_{k,p,q}
 (f)=\sup_{(x,t)\in\Gamma}\left\{(1+x^2+t^2)^k\left|\frac{\partial^{p+q}}{\partial x^p\ \partial t^q}f(x,t)\right|\right\}<+\infty.$$\\
 $\bullet$ $L_\alpha^p(\mathbb{K})$, $p\in[1,+\infty]$, the spaces of measurable functions on $\mathbb{K}$ such that
 \begin{align*}
 &\|f\|_{\alpha,p}=\left[\displaystyle\int_{\mathbb{K}}|f(x,t)|^p\,dm_\alpha(x,t)
\right]^{\frac{1}{p}}<+\infty\,,\ \ \hbox{if}\ p\in[1,+\infty)\\
&\|f\|_{\alpha,\infty}=\mathop{\rm
ess\,sup}_{(x,t)\in \mathbb{K}}|f(x,t)|<+\infty,
 \end{align*}
 where $m_\alpha$ is the positive measure defined on $\mathbb{K}$ by
$$dm_\alpha(x,t)=\frac{1}{\pi\Gamma(\alpha+1)}x^{2\alpha+1}\ t^{2\alpha}\ dx\ dt.$$
 $\bullet$ $L_{\gamma_\alpha}^p(\widehat{\mathbb{K}})$, $p\in[1,+\infty]$, the spaces of measurable functions on $\widehat{\mathbb{K}}$ such that
 \begin{align*}
 &\|g\|_{\gamma_\alpha,p}=\left[\displaystyle\int_{\widehat{\mathbb{K}}}|g(\lambda,m)|^p\,d\gamma_\alpha(\lambda,m)
\right]^{\frac{1}{p}}<+\infty\,,\ \ \hbox{if}\ p\in[1,+\infty)\\
&\|g\|_{\gamma_\alpha,\infty}=\mathop{\rm
ess\,sup}_{(\lambda,m)\in \widehat{\mathbb{K}}}|g(\lambda,m)|<+\infty,
 \end{align*}
 where $\gamma_\alpha$ is the positive measure defined on $\widehat{\mathbb{K}}$ by
$$\int_{[0,+\infty[\times \mathbb{N}}g(\lambda,m)\ d\gamma_\alpha(\lambda,m)=\frac{1}{2^{2\alpha-1}\Gamma(\alpha+\frac{1}{2})}
\sum_{m=0}^{\infty}L_m^\alpha(0)\displaystyle\int_{0}^{+\infty}g(\lambda,m)\ \ \lambda^{3 \alpha+1}d\lambda.$$
Let $f\in \mathcal{S}_*(\mathbb{K})$, for all $(x,t)$ and $(y,s)\in\mathbb{K}$, we put\\
$$T_{(x,t)}^{(\alpha)}f(y,s)=\left\{\begin{array}{lll}\frac{1}{4\pi}\displaystyle\sum_{i,j=0}^{1}\int_{0}^{\pi}f(\Delta_\theta(x,y),Y+(-1)^it+(-1)^js)\ d\theta\ ,\qquad \\\\
 if \ \alpha=0,\,\\  \\ b_\alpha \ \int_{[0,\pi]^3}f\left(\Delta_\theta(x,y),\Delta_\theta(x,y)\xi\right)\ \ d\mu_\alpha(,\xi,\psi,\theta)\ ,\\\\
 if \ \alpha > 0.\end{array}\right.$$
where   $\Delta_\theta(x,y)=\sqrt{x^2+y^2+2xy\cos\theta}$, $b_\alpha=\displaystyle\frac{(\alpha+1)\Gamma(\alpha+\frac{1}{2})}{\Pi^3\Gamma(\alpha)}$,\\ $Y=xy\sin\theta$ and $$dm_\alpha(,\xi,\psi,\theta)=(\sin\xi)^{2\alpha-1}\ (\sin\psi)^{2\alpha-1}\ (\sin\theta)^{2\alpha}d\xi\ d\psi\ d\theta.$$
We define the convolution priduct $f\ast g$ of two functions $f,g\in\mathcal{S}_*(\mathbb{K})$, by
$$(f\ast g)(x,t)=\int_{\mathbb{K}}T_{(x,t)}^{(\alpha)}f(y,s)\ g(y,s)dm_\alpha(y,s),\ \ \ (x,t)\in \mathbb{K}.$$
\begin{lem}
If $f\in L_\alpha^p(\mathbb{K})$, $g\in L_\alpha^q(\mathbb{K})$ such that $1\leq p,q \leq \infty$ and \\ $\frac{1}{p}+\frac{1}{q}-1=\frac{1}{r}$, then the function $f\ast g\ \in L_\alpha^r(\mathbb{K})$, and
$$\left\|f \ast g\right\|_{\alpha,r}\leq\ \left\|f \right\|_{\alpha,p}\ \left\|g \right\|_{\alpha,q}.$$
\end{lem}
We consider the dilations on $\mathbb{K}$ defined by
$$\delta_r(x,t)=(rx,r^2t),\ \ \ \ \ \ \ \ \ \ \ r>0.$$
We also introduce a homogeneous norm, related to family $(\delta_r)_{r>0}$ defined by $$\left|(x,t)\right|=(x^4+4t^2)^{\frac{1}{4}}$$
We define the ball centered at $(0,0)$ of radius $r$ by  $$B_r=\left\{(x,t) \in \mathbb{K} ; \left|(x,t)\right|< r \right\}.$$
Let $f\in L_\alpha^1(\mathbb{K})$, the Laguerre-Bessel transform of $f$ is defined by $$\mathcal{F}_{LB}(f)(\lambda,m)=\displaystyle\int_{\mathbb{K}}f(x,t)\ \varphi_{(\lambda,m)}(x,t)dm_\alpha(x,t).$$
For $f$ and $g\in L_\alpha^1(\mathbb{K})$, we have :
$$\mathcal{F}_{LB}(f\ast g)(\lambda,m)=\mathcal{F}_{LB}(f)(\lambda,m) \ \mathcal{F}_{LB}(g)(\lambda,m).$$
The integral transform can be extended to an isometric isomorphism
$L_\alpha^2(\mathbb{K})$ to $L_{\gamma_\alpha}^2(\widehat{\mathbb{K}})$ and
we have the Plancherel formula.$$\left\|f\right\|_{\alpha,2}=\left\|\mathcal{F}_{LB}(f)\right\|_{\gamma_\alpha,2},\ \ \ \ \ f\in L^1(\mathbb{K})\cap L^2(\mathbb{K}).$$
We consider the differential operator
\begin{equation}\label{lh}
L=-\left(\frac{\partial^2}{\partial x^2}+\frac{2\alpha+1}{x}\frac{\partial}{\partial x}+x^2D_1\right).
\end{equation}
$L$ is positive and symetric in $L_\alpha^2(\mathbb{K})$, and is homogeneous of degree $2$ if $\mathbb{K}$ is endowed with the family of dilations $(\delta_r)_{r>0}$, $\delta_r(x,t)=(rx,r^2t)$.\\
We have $$L\varphi_{(\lambda,m)}=2\lambda(2m+\alpha+1)\varphi_{(\lambda,m)}.$$
As in \cite{Ste}, page $117$, we define $L^{b}$ for $b>0$ by $$\mathcal{F}_{LB}\left(L^{b}f\right)(\lambda,m)=\left(2 \lambda(2m+\alpha+1)\right)^{b}\ \ \mathcal{F}_{LB}(f)(\lambda,m),$$
On the other hand,
$L$ is hypoelliptic on $\mathbb{K}$. Also, the heat operator\\ $L+\partial_s$ is hypoelliptic on $\mathbb{K}\times (0,+\infty)$. Hence, similar arguments from the proof of Hunt's theorem \cite[Theorem 3.4]{Hun}
\begin{prop}
There is a unique $C^\infty$ function $h\left((x,t),s\right)=h_s(x,t)$ on $\mathbb{K}\times(0,+\infty)$ with the following properties\\
    i) $(L+\partial_s)h=0$ on $\mathbb{K}\times(0,+\infty)$,\\
    ii) $h_s(x,t)\geq0$ and $\displaystyle\int_{\mathbb{K}}h_s \ dm_\alpha=1$,\\
    iii) $h_{s_1}\ast h_{s_2}=h_{s_1+s_2}$,       $s_1,s_2>0$,\\
\end{prop}
\begin{lem}
For any $s>0$, $\mathcal{F}_{LB}(h_s)(\lambda,m)=e^{-2\lambda(2m+\alpha+1)s}$.
\end{lem}
\textbf{Proof}.
From the equalities
\begin{equation}\label{ll}
\partial_s (h_s\ast u)=-Lu\ast h_s
\end{equation}
and
\begin{equation}\label{lll}
\mathcal{F}_{LB}(Lu)(\lambda,m)=-2\lambda(2m+\alpha+1) \ \mathcal{F}_{LB}(u)(\lambda,m).
\end{equation}
we show that the function $\mathcal{F}_{LB}(h_s)(\lambda,m)$ satisfy the differential equations $\displaystyle\frac{d}{ds}w=-2\lambda(2m+\alpha+1)w.$
the result is proved.\hfill$\square$\\
Let $\left\{H^s, \ \  s >0\right\}$ be the heat semigroup. There is an unique smooth function $h((x,t),s)=h_s(x,t)$ on $\mathbb{K}\times (0,+\infty)$ such that
$H^sf(x,t)=f\ast h_s(x,t)$.\\
$h_s$ is called the heat kernel assocaited to $L$.
\begin{lem}
\begin{equation}\label{heisen}
\left\|h_s\right\|_{\alpha,2}\leq C\ s^{-\frac{3\alpha+2}{2}}
\end{equation}
\end{lem}
\textbf{Proof}.
By the Plancherel formula, we have $\left\|h_s\right\|_{\alpha,2}=\left\|\mathcal{F}_{LB}(h_s)\right\|_{\gamma_\alpha,2}$.
$$\left\|\mathcal{F}_{LB}(h_s)\right\|_{\gamma_\alpha,2}^2=\displaystyle\frac{1}{2^{2\alpha-1}\Gamma(\alpha+\frac{1}{2})}\displaystyle\int_{0}^{+\infty}\left(\displaystyle\sum_{m=0}^{+\infty}L_m^\alpha(0)e^{-8\lambda s m}\right)\ e^{-4\lambda s (\alpha+1)} \ \lambda ^{3\alpha+1} \ d\lambda\\
$$
By the generating function identity $(\ref{lag})$ for the Laguerre polynomials, we have:
$$
\left\|\mathcal{F}_{LB}(h_s)\right\|_{\gamma_\alpha,2}^2=\displaystyle\frac{s^{-3\alpha-2}}{2^{2\alpha-1}\Gamma(\alpha+\frac{1}{2})}\displaystyle\int_{0}^{+\infty}\left(\displaystyle\frac{1}{2\sinh(4u)}\right)^{\alpha+1}\  \ u ^{3\alpha+1}  \ du\\
$$
So,\qquad $\left\|\mathcal{F}_{LB}(h_s)\right\|_{\gamma_\alpha,2}^2\leq C \ s^{-(3\alpha+2)}\,.$\hfill$\square$
\section{Heisenberg inequality for Laguerre-Bessel transform}
\begin{lem}\label{heilag}
Let $0<a< 3\alpha+2$, then for all $f\in L_\alpha^2(\mathbb{K})$, we have $$\left\|H^sf\right\|_{\alpha,2}\leq \ C \ \ s^{-\frac{a}{2}}\ \ \left\|\ |(x,t)|^a\ f\right\|_{\alpha,2}.$$
\end{lem}
\textbf{Proof}.
For $r>0$, let $f_r=f\ \chi_{B_r}$ and $f^r=f-f_r$.\\
Then $$\left|f^r(x,t)\right|\leq \ r^{-a}\  |(x,t)|^a\ \left|\ f(x,t)\right|$$
So
$$\left\|H^sf^r\right\|_{\alpha,2}\leq \ \left\| f^r\right\|_{\alpha,2}\leq \ r^{-a}  \left\|\ |(x,t)|^a\ f\right\|_{\alpha,2}.$$
On the other hand, we have
\begin{align*}
\left\|H^sf_r\right\|_{\alpha,2}&=\ \left\| f_r\ast h_s\right\|_{\alpha,2}\\
&\leq \ \  \left\|f_r\right\|_{\alpha,1}\  \left\|h_s\right\|_{\alpha,2}\\
&\leq \ \left\|h_s\right\|_{\alpha,2} \ \left\|\ |(x,t)|^{-a}\chi_{B_r}\right\|_{\alpha,2}\ \  \left\|\ |(x,t)|^a\ f\right\|_{\alpha,2}
\end{align*}
Since,
$$\left\|\ |(x,t)|^{-a}\chi_{B_r}\right\|^2_{\alpha,2}=\displaystyle\frac{B(\frac{\alpha+1}{2},\frac{2\alpha+1}{2})}{4^{\alpha+1}\pi
\Gamma(\alpha+1)(3\alpha+2-a)}\  \ r^{6\alpha+4-2a}$$
with $B$ is the beta function, we get
\begin{align*}
\left\|H^sf\right\|_{\alpha,2} & \leq \left\|H^sf_r\right\|_{\alpha,2}+\left\|H^sf^r\right\|_{\alpha,2}\\
&\leq \ r^{-a}  \left\| \ |(x,t)|^a\ f\right\|_{\alpha,2} \left(1+C\ \left\|h_s\right\|_{\alpha,2}\ r^{3\alpha+2}\right)
\end{align*}
By the relation (\ref{heisen}), we obtain:
$$\left\|H^sf\right\|_{\alpha,2}
\leq \ r^{-a}  \left\| \ |(x,t)|^a\ f\right\|_{\alpha,2} \left(1+C\ s^{-\frac{3\alpha+2}{2}}\ r^{3\alpha+2}\right)$$
Choosing $r=s^{\frac{1}{2}}$,
we obtain\quad$\left\|H^sf\right\|_{\alpha,2}\leq \ C\ s^{-\frac{a}{2}}  \left\|\  |(x,t)|^a\ f\right\|_{\alpha,2}\,.$\hfill$\square$
\begin{thm}
Let $a,b>0$, then for all $f\in L_\alpha^2(\mathbb{K})$, we have
\begin{equation}\label{th}
\left\|\  |(x,t)|^a\ f\right\|_{\alpha,2}^{\frac{2b}{a+2b}} \ \ \left\|\left(2(2m+\alpha+1)\lambda\right)^{b}\mathcal{F}_{LB}(f)
\right\|_{\gamma_\alpha,2}^{\frac{a}{a+2b}} \geq C\ \left\|f\right\|_{\alpha,2}
\end{equation}
\end{thm}
\textbf{Proof }.
$\mathcal{S}(\mathbb{K})$ is dense in $L_\alpha^2(\mathbb{K})$, so we need only to prove $(\ref{th})$ for $\mathcal{S}(\mathbb{K}).$ \\
Assume that $a<3\alpha+2$.\\
If $b\leq 1$,
By lemma \ref{heilag},
\begin{align*}
\left\| f\right\|_{\alpha,2}\leq&\left\|H^{s}f\right\|_{\alpha,2}+\left\|(1-H^{s})\ f\right\|_{\alpha,2}\\
&\leq C \ s^{-\frac{a}{2}}  \left\|\ |(x,t)|^a \ f\right\|_{\alpha,2} + \left\|(1-H^{s})\ \ (s L)^{-b}\ (s L)^{b}f\right\|_{\alpha,2}
\end{align*}
Let $g=(s L)^{b}f$, so
$$
\left\|(1-H^{s})\ \ (s L)^{-b}\ g \right\|_{\alpha,2}
=\left\|\ (1-e^{2\lambda(2m+\alpha+1)})\ (2\lambda(2m+\alpha+1))^{-b} \ \ \mathcal{F}_{LB}\left(g\right)\right\|_{\gamma_\alpha,2}
$$
since, if $b\leq1$ the function $t\longmapsto(1-e^{-t}) \ t^{-b}$ is bounded for $t\geq 0$ . Therefore $$\left\| f\right\|_{\alpha,2}
\leq C \left(\ s^{-\frac{a}{2}}  \left\|\ |(x,t)|^a \ f\right\|_{\alpha,2} + s^{
b}\left\|L^{b}\ f\right\|_{\alpha,2}\right)$$
From which, optimizing in $s$, we obtain $$\left\|\ |(x,t)|^a \ f\right\|_{\alpha,2}^{\frac{2b}{a+2b}} \ \ \left\|L^{b}f
\right\|_{\alpha,2}^{\frac{a}{a+2b}}\geq C\ \left\|f\right\|_{\alpha,2}$$
Since $\mathcal{F}_{LB}(L^{b}f)(\lambda,m)=\left(2(2m+\alpha+1)\lambda\right)^{b}\mathcal{F}_{LB}(f)(\lambda,m)$ and from the Plancherel formula.
we get the result.\\
If $b>1$. For $u\geq0$, $u\leq 1+u^b$, which for
$u=\frac{2(2m+\alpha+1)\ \lambda}{\varepsilon}$ gives the \\\\
inequality $\frac{2(2m+\alpha+1) \lambda}{\varepsilon}\leq 1+ (\frac{2(2m+\alpha+1)\lambda}{\varepsilon})^{b}$, for all $\varepsilon>0$.\\\\
It follows that $$\left\|\left(2(2m+\alpha+1)\lambda\right)\mathcal{F}_{LB}(f)\right\|_{\gamma_\alpha,2}\leq\ \varepsilon\left\|f\right\|_{\alpha,2}+ \varepsilon^{1-b}\left\|\left(2(2m+\alpha+1)\lambda\right)^{b}\mathcal{F}_{LB}(f)\right\|_{\gamma_\alpha,2}
$$
optimizing in $\varepsilon$, we get:
$$\left\|\left(2(2m+\alpha+1)\lambda\right)\mathcal{F}_{LB}(f)\right\|_{\gamma_\alpha,2}\leq\ C \left\|f\right\|_{\alpha,2}^{1-\frac{1}{b}}\left\|\left(2(2m+\alpha+1)\lambda\right)^{b}\mathcal{F}_L(f)\right\|_{\gamma_\alpha,2}^{\frac{1}{b}}
$$
Together with $(\ref{th})$ for $b=1$, we get the result for $b> 1$.\\
If $a\geq 3\alpha+2$, then using $$
\displaystyle\frac{|(x,t)|}{\varepsilon}\leq 1+\displaystyle\frac{|(x,t)|^a}{\varepsilon^{a}},\ \ \ \varepsilon>0,$$
It follows that $$\left\| \
|(x,t)| \ f\right\|_{\alpha,2}\leq\ \varepsilon\left\|f\right\|_{\alpha,2}+ \varepsilon^{1-a}\left\|\ |(x,t)|^a \  f\right\|_{\alpha,2}
$$
optimizing in $\varepsilon$, we get:
$$\left\|\ |(x,t)| \ f\right\|_{\alpha,2}\leq\ C \left\|f\right\|_{\alpha,2}^{1-\frac{1}{a}}\left\|\ |(x,t)|^a \ f\right\|_{\alpha,2}^{\frac{1}{a}}
$$
Together with $(\ref{th})$ for $a=1$, we get the result for $a\geq 3\alpha+2$.\hfill$\square$
\section{Local uncertainty inequalities}
In this section, we establish a local uncertainty inequalities related to Laguerre-Bessel transform.
Similar results are obtained by Omri and Rachdi $\cite{omr}$ in framework of the
Riemann-Liouville operator.
\begin{thm}
Let $s$ be a real number such that $0<s< 3\alpha+2$. Then for all nonzero $f\in L_\alpha^2(\mathbb{K})$ and for all measurable subsets $E\subset \widehat{\mathbb{K}}$ such that $0<\gamma_\alpha(E)<+\infty$, we have
\begin{equation}\label{1.1}
\left(\displaystyle\int_{}\displaystyle\int_{E}|\mathcal{F}_{LB}(f)(\lambda,m)|^2d\gamma_\alpha(\lambda,m)\right)^\frac{1}{2}< K_{\alpha,s}\ \  \gamma_\alpha(E)^{\frac{s}{2(3\alpha+2)}}\ \left\|\ |(x,t)|^s \ f\right\|_{\alpha,2}
\end{equation}
    where $$K_{\alpha,s}=\left(\frac{B(\frac{\alpha+1}{2},\frac{2\alpha+1}{2})(3\alpha+2-s)}{4^{\alpha+1}\pi\Gamma(\alpha+1)s^2}\right)^{\frac{s}{2(3\alpha+2)}}
    \displaystyle\frac{3\alpha+2}{3\alpha+2-s}$$
\end{thm}
\textbf{Proof}.
Let $0<s< 3\alpha+2$ and $f\in L_\alpha^2(\mathbb{K})$, we have
$$\left(\displaystyle\int_{}\displaystyle\int_{E}|\mathcal{F}_{LB}(f)(\lambda,m)|^2d\gamma_\alpha(\lambda,m)\right)^{\frac{1}{2}}=\left\|\mathcal{F}_{LB}(f)\chi_E\right\|_{\gamma_\alpha,2}$$
 By Minkowski's inequality, it follows $$\left\|\mathcal{F}_{LB}(f)\chi_E\right\|_{\gamma_\alpha,2}\leq \left\|\mathcal{F}_{LB}(f\chi_{B_r})\chi_E\right\|_{\gamma_\alpha,2}+\left\|\mathcal{F}_{LB}(f\chi_{B_r^c})\chi_E\right\|_{\gamma_\alpha,2}$$
 Therefore
  \begin{equation}\label{1.2}
\left\|\mathcal{F}_{LB}(f)\chi_E\right\|_{\gamma_\alpha,2} \leq \gamma_\alpha(E)^{\frac{1}{2}} \left\|\mathcal{F}_{LB}(f\chi_{B_r})\right\|_{\gamma_\alpha,\infty}+\left\|\mathcal{F}_{LB}(f\chi_{B_r^c})\right\|_{\gamma_\alpha,2}
\end{equation}
\begin{equation}\label{1.3}
\left\|\mathcal{F}_{LB}(f)\chi_E\right\|_{\gamma_\alpha,2} \leq \gamma_\alpha(E)^{\frac{1}{2}} \left\|f\chi_{B_r}\right\|_{\alpha,1}+\left\|\mathcal{F}_{LB}(f\chi_{B_r^c})\right\|_{\gamma_\alpha,2}
\end{equation}
On the other hand, using H$\ddot{o}$lder inequality, we get
\begin{equation}\label{1.4}
\left\|f\chi_{B_r}\right\|_{\alpha,1}\leq
\left\| \ |(x,t)|^s \ f\right\|_{\alpha,2}\ \left\| \ |(x,t)|^{-s}\chi_{B_r}\right\|_{\alpha,2}
\end{equation}
Therefore, we have
\begin{equation}\label{1.5}
\left\|f\chi_{B_r}\right\|_{\alpha,1}\leq \left\| \ |(x,t)|^s \ f\right\|_{\alpha,2}\ \left(\displaystyle\frac{B(\frac{\alpha+1}{2},\frac{2\alpha+1}{2})}{4^{\alpha+1} \pi\Gamma(\alpha+1)(3\alpha+2-s)}\right)^{\frac{1}{2}}r^{3\alpha+2-s}
\end{equation}
Plancherel's theorem allows as to say
\begin{align*}
\left\|\mathcal{F}_{LB}(f\chi_{B_r^c})\right\|_{\gamma_\alpha,2}&=\left\|f\chi_{B_r^c}\right\|_{\alpha,2}\\
&\leq \left\| \ |(x,t)|^s \ f\right\|_{\alpha,2}\ \left\| \ |(x,t)|^{-s} \ \chi_{B_r^c}\right\|_{\alpha,\infty}
\end{align*}
So
\begin{equation}\label{1.6}
\left\|\mathcal{F}_{LB}(f\chi_{B_r^c})\right\|_{\gamma_\alpha,2} \leq r^{-s}\left\| \ |(x,t)|^s \ f\right\|_{\alpha,2}
\end{equation}
Combining the relations (\ref{1.3}),(\ref{1.5}) and (\ref{1.6}), we deduce that for all $s>0$, we have
\begin{equation}\label{1.7}
\left\|\mathcal{F}_{LB}(f)\chi_E\right\|_{\gamma_\alpha,2} \leq g_{\alpha,s}(r)\ \ \left\| \ |(x,t)|^s \ f\right\|_{\alpha,2}
\end{equation}
where $g_{\alpha,s}$is the function defined on $(0,+\infty)$ by $$g_{\alpha,s}(r)=r^{-s}+\left(\displaystyle\frac{B(\frac{\alpha+1}{2},\frac{2\alpha+1}{2})}{4^{\alpha+1}\pi\Gamma(\alpha+1)(3\alpha+2-s)}\gamma_\alpha(E)\right)^{\frac{1}{2}}\ r^{3\alpha+2-s}$$
In particular, we have the inequality (\ref{gs})
\begin{equation}\label{gs}
\left\|\mathcal{F}_{LB}(f)\chi_E\right\|_{\gamma_\alpha,2} \leq g_{\alpha,s}(r_0)\ \ \left\| \ |(x,t)|^s \ f\right\|_{\alpha,2}
\end{equation}
where $$r_0=(\frac{s}{3\alpha+2-s})^{\frac{1}{3\alpha+2}} \ \left(\displaystyle\frac{B(\frac{\alpha+1}{2},\frac{2\alpha+1}{2})}{4^{\alpha+1}\pi\Gamma(\alpha+1)(3\alpha+2-s)}\gamma_\alpha(E)\right)^{-\frac{1}{2(3\alpha+2)}}$$
However $g_{\alpha,s}(r_0)=\  \gamma_\alpha(E)^{\frac{s}{2(3\alpha+2)}}\ K_{\alpha,s}\ $
where $$K_{\alpha,s}=\left(\frac{B(\frac{\alpha+1}{2},\frac{2\alpha+1}{2})\ (3\alpha+2-s)}{4^{\alpha+1}\pi\Gamma(\alpha+1)s^2}\right)^{\frac{s}{2(3\alpha+2)}}
    \displaystyle\frac{3\alpha+2}{3\alpha+2-s}$$
Let us prove that the equality in (\ref{gs}) cannot hold. Suppose that there exists a nonzero function $f\in L_\alpha^2(\mathbb{K})$ such that $$
\left\|\mathcal{F}_{LB}(f)\chi_E\right\|_{\gamma_\alpha,2}= K_{\alpha,s}\ \  \gamma_\alpha(E)^{\frac{s}{2(3\alpha+2)}}\ \left\| \ |(x,t)|^s \ f\right\|_{\alpha,2}$$
Let $$\psi(r)=\gamma_\alpha(E)^{\frac{1}{2}} \left\|\mathcal{F}_{LB}(f\chi_{B_r})\right\|_{\gamma_\alpha,\infty}+\left\|\mathcal{F}_{LB}(f\chi_{B_r^c})\right\|_{\gamma_\alpha,2},\quad r>0$$
We have $$\forall r>0, \quad \psi(r)\leq g_{\alpha,s}(r)\ \ \left\| \ |(x,t)|^s \ f\right\|_{\alpha,2}$$
In particular, $$\psi(r_0)\leq g_{\alpha,s}(r_0)\ \ \left\| \ |(x,t)|^s \ f\right\|_{\alpha,2}$$
But
\begin{equation}\label{1.111}
\left\|\mathcal{F}_{LB}(f)\chi_E\right\|_{\gamma_\alpha,2}=g_{\alpha,s}(r_0)\ \ \left\| \ |(x,t)|^s \ f\right\|_{\alpha,2}
\end{equation}
So by the relations (\ref{1.2}) and (\ref{1.111}), we get
\begin{equation}\label{1.122}
\left\|\mathcal{F}_{LB}(f)\chi_E\right\|_{\gamma_\alpha,2}=\gamma_\alpha(E)^{\frac{1}{2}}\ \ \left\|\mathcal{F}_{LB}(f\chi_{B_{r_0}})\right\|_{\gamma_\alpha,\infty}+\left\|\mathcal{F}_{LB}(f\chi_{B_{r_0}^c})\right\|_{\gamma_\alpha,2}
\end{equation}
On the other hand, we have
\begin{equation}\label{1.113}
\left\|\mathcal{F}_{LB}(f)\chi_E\right\|_{\gamma_\alpha,2}\leq \ \left\|\mathcal{F}_{LB}(f\chi_{B_{r_0}})\chi_E\right\|_{\gamma_\alpha,2}+\left\|\mathcal{F}_{LB}(f\chi_{B_{r_0}^c})\chi_E\right\|_{\gamma_\alpha,2}
\end{equation}
Using the relations (\ref{1.122}) and (\ref{1.113}), we have
\begin{equation}\label{1.114}
\gamma_\alpha(E)^{\frac{1}{2}}\ \ \left\|\mathcal{F}_{LB}(f\chi_{B_{r_0}})\right\|_{\gamma_\alpha,\infty}\leq
\left\|\mathcal{F}_{LB}(f\chi_{B_{r_0}})\chi_E\right\|_{\gamma_\alpha,2}
\end{equation}
Writting the relation (\ref{1.2}) for the function $f\chi_{B_{r_0}}$, we obtain
\begin{equation}\label{1.115}
 \left\|\mathcal{F}_{LB}(f\chi_{B_{r_0}})\chi_E\right\|_{\gamma_\alpha,2}\leq \gamma_\alpha(E)^{\frac{1}{2}}\ \ \left\|\mathcal{F}_{LB}(f\chi_{B_{r_0}})\right\|_{\gamma_\alpha,\infty},
\end{equation}
Therefore
\begin{equation}\label{1.8}
 \left\|\mathcal{F}_{LB}(f\chi_{B_{r_0}})\chi_E\right\|_{\gamma_\alpha,2}= \gamma_\alpha(E)^{\frac{1}{2}}\ \ \left\|\mathcal{F}_{LB}(f\chi_{B_{r_0}})\right\|_{\gamma_\alpha,\infty}.
\end{equation}
Combining the relations (\ref{1.3}),(\ref{1.5}) and (\ref{1.6}), we obtain
\begin{equation}\label{rr}
\gamma_\alpha(E)^{\frac{1}{2}} \left\|f\chi_{B_r}\right\|_{\alpha,1}+\left\|\mathcal{F}_L(f\chi_{B_r^c})\right\|_{\gamma_\alpha,2}\leq g_{\alpha,s}(r)\ \ \left\| \ |(x,t)|^s \ f\right\|_{\alpha,2};\ r>0
\end{equation}
The relations (\ref{1.3}), (\ref{1.111}) and (\ref{rr}), lead
\begin{equation}\label{1.110}
 \left\|\mathcal{F}_{LB}(f)\chi_E\right\|_{\gamma_\alpha,2}=\gamma_\alpha(E)^{\frac{1}{2}} \left\|f\chi_{B_{r_0}}\right\|_{\alpha,1}+\left\|\mathcal{F}_{LB}(f\chi_{B_{r_0}^c})\right\|_{\gamma_\alpha,2}
\end{equation}
So, using (\ref{1.122}) we get
\begin{equation}\label{1.9}
\left\|f\chi_{B_{r_0}}\right\|_{\alpha,1}=
\left\|\mathcal{F}_{LB}(f\chi_{B_{r_0}})\right\|_{\gamma_\alpha,\infty}
\end{equation}
Using the relations (\ref{1.3}) and (\ref{1.4}), we have
\begin{equation}\label{1.124}
 \left\|\mathcal{F}_{LB}(f)\chi_E\right\|_{\gamma_\alpha,2}\leq \varphi(r);\quad r>0
 \end{equation}
with
\begin{equation}\label{rrr}
\varphi(r)=\gamma_\alpha(E)^{\frac{1}{2}}\left\| \ |(x,t)|^s \ f\right\|_{\alpha,2} \left\| \ |(x,t)|^{-s} \ \chi_{B_{r}}\right\|_{\alpha,2}+\left\|\mathcal{F}_{LB}(f\chi_{B_r^c})\right\|_{\gamma_\alpha,2}
\end{equation}
We have $$\forall r>0, \quad \varphi(r)\leq g_{\alpha,s}(r)\ \ \left\| \ |(x,t)|^s \ f\right\|_{\alpha,2}$$
In particular , $$\varphi(r_0)\leq g_{\alpha,s}(r_0)\ \ \left\| \ |(x,t)|^s \ f\right\|_{\alpha,2}$$
But
$$g_{\alpha,s}(r_0)=\displaystyle\frac{\left\|\mathcal{F}_{LB}(f)\chi_E\right\|_{\gamma_\alpha,2}}{ \left\| \ |(x,t)|^s \ f\right\|_{\alpha,2}}$$ Therefore
\begin{equation}\label{1.115}
\varphi(r_0)\leq \left\|\mathcal{F}_{LB}(f)\chi_E\right\|_{\gamma_\alpha,2}
\end{equation}
Using the relations (\ref{1.110}), (\ref{1.124}), (\ref{rrr}) and (\ref{1.115}), we have
\begin{equation}\label{1.10}
\left\|f\chi_{B_{r_0}}\right\|_{\alpha,1}=
\left\| \ |(x,t)|^s \ f\right\|_{\alpha,2}\ \left\| \ |(x,t)|^{-s} \ \chi_{B_{r_0}}\right\|_{\alpha,2}
\end{equation}
However, $f$ satisfies the equality (\ref{1.10}) if and only if
$$ \left|f(x,t)\right|=C\  |(x,t)|^{-2s}\chi_{B_{r_0}}(x,t),$$
hence
\begin{equation}\label{1.11}
\forall (x,t)\in \mathbb{K}, \ f(x,t)= C \ \phi(x,t) \ |(x,t)|^{-2s}\chi_{B_{r_0}}(x,t),
\end{equation}
with $|\phi(x,t)|=1$.\\
But $f$ satisfies the relation (\ref{1.9}), then there exists $(\lambda_0,m_0)\in \hat{\mathbb{K}}$, such that
$$ \left\|f\right\|_{\alpha,1}=\left\|\mathcal{F}_{LB}(f)\right\|_{\gamma_\alpha,\infty}=\left|\mathcal{F}_{LB}(f)(\lambda_0,m_0)\right|.$$
So, there exists $\theta_0 \in \mathbb{R}$ satisfying
\begin{equation}\label{1.80}
\mathcal{F}_{LB}(f)(\lambda_0,m_0)=e^{i\theta_0}\left\|f\right\|_{\alpha,1},
\end{equation}
and therefore
$$Ce^{i\theta_0}\displaystyle\int_{\mathbb{K}} \ |(x,t)|^{-2s}\chi_{B_{r_0}}(x,t)\left(\ \phi(x,t)e^{-i\theta_0}\ j_{\alpha-\frac{1}{2}}(\lambda_0t)\  \mathcal{L}_{m_0}^\alpha(\lambda_0x^2)-1\right)dm_\alpha(x,t)=0.$$
This implies that for almost every $(x,t)\in \mathbb{K}$,
$$\phi(x,t)\ e^{-i\theta_0}\ j_{\alpha-\frac{1}{2}}(\lambda_0t)\  \mathcal{L}_{m_0}^\alpha(\lambda_0 x^2)=1.$$
Since $|\phi(x,t)|=1$, we deduce that for all $x\in \mathbb{R}_+$,
$$\left|\ j_{\alpha-\frac{1}{2}}(\lambda_0t)\ \mathcal{L}_{m_0}^\alpha(\lambda_0 x^2)\right|=1.$$
It follows that $\lambda_0=0$ and then $$\phi(x,t)=e^{i\theta_0}.$$
Replacing in (\ref{1.11}), we get $$f(x,t)=C\ e^{i\theta_0}\  \ |(x,t)|^{-2s}\chi_{B_{r_0}}(x,t).$$
On the other hand, by the relation (\ref{1.8}), we get
$$\displaystyle\int_{}\displaystyle\int_{E}(\left\|\mathcal{F}_{LB}(f)\right\|^2_{\gamma_\alpha,\infty}-\left|\mathcal{F}_{LB}(f)(\lambda,m)|^2\right)d\gamma_\alpha(\lambda,m)=0$$
 then for almost every $(\lambda,m)\in E$, we have
$$\left|\mathcal{F}_{LB}(f)(\lambda,m)\right|=\left\|\mathcal{F}_{LB}(f)\right\|_{\gamma_\alpha,\infty},$$
and by (\ref{1.80}), we deduce that for almost every $(\lambda,m)\in E$,
$$\left|\mathcal{F}_{LB}(f)(\lambda,m)\right|=e^{-i\theta_0}\ \mathcal{F}_{LB}(f)(0,m_0)$$
Hence,
$$\mathcal{F}_{LB}(f)(\lambda,m)=\varphi(\lambda,m)\ e^{-i\theta_0}\ \mathcal{F}_{LB}(f)(0,m_0),$$
with $|\varphi(\lambda,m)|=1$, and therefore $$C\displaystyle\int_{\mathbb{K}} \ |(x,t)|^{-2s}\chi_{B_{r_0}}(x,t)\left(\varphi(\lambda,m)^{-1}\ e^{i\theta_0}\ \ j_{\alpha-\frac{1}{2}}(\lambda t)\ \mathcal{L}_m^\alpha(\lambda x^2)-1\right)dm_\alpha(x,t)=0.$$
Consequently for almost every
$(x,t)\in \mathbb{K}$,
$$\varphi(\lambda,m)^{-1}\ e^{i\theta_0}\ \ j_{\alpha-\frac{1}{2}}(\lambda t)\ \mathcal{L}_m^\alpha(\lambda x^2)=1.$$
which implies that $\lambda=0$.
However, since $\gamma_\alpha(E)>0$, this contradicts the fact that for almost every $(\lambda,m)\in E$,
$$\left|\mathcal{F}_{LB}(f)(\lambda,m)\right|=\left|\mathcal{F}_{LB}(f)(0,m_0)\right|,$$
and shows that the inequality in (\ref{1.1}) is stictly satisfied.\hfill$\square$
\begin{lem}\label{5.2}
Let $s$ be a real number such that $s>3\alpha+2$, then for all nonzero measurable function $f$ on $\mathbb{K}$, we have
\begin{equation}\label{1.12}
\left\|f\right\|_{\alpha,1}\leq M_{\alpha,s}\ \left\|f\right\|^{1-\frac{3\alpha+2}{s}}_{\alpha,2}\ \left\| \ |(x,t)|^s \ f\right\|^{\frac{3\alpha+2}{s}}_{\alpha,2},
\end{equation}
where  $M_{\alpha,s}=\left(\frac{B(\frac{\alpha+1}{2},\frac{2\alpha+1}{2})\ B(\frac{s-3\alpha-2}{s},\frac{3\alpha+2}{s})}{4^{\alpha+1}\pi\Gamma(\alpha+1)(s-3\alpha-2)}\ \left(\frac{s-3\alpha-2}{3\alpha+2}\right)^{\frac{3\alpha+2}{s}}\right)^{\frac{1}{2}}$.\\
We have equality in $(\ref{1.12})$ if only if there exists $a>0$ and $b>0$ such that: $$|f(x,t)|=\left(a+b \ |(x,t)|^{2s} \right)^{-1}$$
\end{lem}
\textbf{Proof}.
The inequality (\ref{1.12}) holds if $\left\|f\right\|_{\alpha,2}=+\infty$ or $\left\| \ |(x,t)|^s \ f\right\|_{\alpha,2}=+\infty$.\\\\
Assume that $\left\|f\right\|_{\alpha,2}+\left\| \ |(x,t)|^s \ f\right\|_{\alpha,2}<+\infty$.\\\\
From the hypothesis $s>3\alpha+2$, we deduce that for all $a>0$ and $b>0$, the function
$$(x,t)\longmapsto \left(a+b\ |(x,t)|^{2s} \right)^{-1}$$
belongs to $L^1_\alpha(\mathbb{K})$ and by H$\ddot{o}$lder's inequality, we have
\begin{equation}\label{1.13}
\left\|f\right\|^2_{\alpha,1}\leq \left\|\left(1+ \ |(x,t)|^{2s}\right)^{\frac{1}{2}}f\right\|^{2}_{\alpha,2}\ \left\|\left(1+ \ |(x,t)|^{2s} \ \right)^{-\frac{1}{2}}\right\|^{2}_{\alpha,2}\
\end{equation}
We have equality in (\ref{1.13}) if and only if
\begin{equation}\label{1.15}
|f(x,t)|=C\  \left(1+ \ |(x,t)|^{2s}\right)^{-1}.
\end{equation}
But $$\left\|\left(1+ \ |(x,t)|^{2s}
\right)^{\frac{1}{2}}f\right\|^{2}_{\alpha,2}=\left\|f\right\|^2_{\alpha,2}+\left\|
\ |(x,t)|^s \ f\right\|^{2}_{\alpha,2} $$ Therefore
\begin{equation}\label{1.14}
\left\|f\right\|^2_{\alpha,1}\leq N_{\alpha ,s} \  \left(\left\|f\right\|^2_{\alpha,2}+\left\| \ |(x,t)|^s \ f\right\|^{2}_{\alpha,2}\right)
\end{equation}
where
$$N_{\alpha,s}=\left\|\left(1+\ |(x,t)|^{2s} \right)^{-\frac{1}{2}}\right\|^{2}_{\alpha,2}.$$
By straightforward calculus, we get $$N_{\alpha,s}=\displaystyle\frac{1}{4^{\alpha+1}\ s \pi\Gamma(\alpha+1)}\ B(\frac{\alpha+1}{2},\frac{2\alpha+1}{2})\
B(\frac{s-3\alpha-2}{s},\frac{3\alpha+2}{s})$$
For $r>0$, we put $$f_r(x,t)=r^{-(6\alpha+4)}f(\frac{x}{r},\frac{t}{r^2}).$$
Then we have
\begin{equation}\label{1.166}
\left\|f_r\right\|_{\alpha,1}=\left\|f\right\|_{\alpha,1}.
\end{equation}
\begin{equation}\label{1.160}
\left\|f_r\right\|^2_{\alpha,2}=\displaystyle\frac{1}{r^{6\alpha+4}}\left\|f\right\|^2_{\alpha,2}
\end{equation}
\begin{equation}\label{1.17}
\left\|\ |(x,t)|^s \ f_r\right\|^2_{\alpha,2}=\displaystyle\frac{1}{r^{6\alpha+4-2s}}\left\|\ |(x,t)|^s \ f\right\|^2_{\alpha,2}
\end{equation}
Replacing $f$ by $f_r$ in the relation (\ref{1.14}), we deduce that for all $r>0$, we have
$$\left\|f\right\|^2_{\alpha,1}\leq N_{\alpha ,s}\  \left(r^{-(6\alpha+4)}\ \left\|f\right\|^2_{\alpha,2}+\ r^{2s-6\alpha-4}\ \left\|\ |(x,t)|^s \
f\right\|^{2}_{\alpha,2}\right).$$
In particular, for
$$r_0=\left(\displaystyle\frac{(3\alpha+2)\left\|f \right\|^2_{\alpha,2}}{(s-3\alpha-2) \left\| \ |(x,t)|^s \ f \right\|^{2}_{\alpha,2}}\right)^{\frac{1}{2s}}$$
we get
\begin{equation}\label{1.227}
\left\|f\right\|^2_{\alpha,1}\leq M^2_{\alpha,s}\ \left\|f\right\|^{2-\frac{6\alpha+4}{s}}_{\alpha,2}\ \left\|\ |(x,t)|^s \  f\right\|^{\frac{6\alpha+4}{s}}_{\alpha,2},
\end{equation}
where \\ $M_{\alpha,s}=\left(\frac{B(\frac{\alpha+1}{2},\frac{2\alpha+1}{2})\ B(\frac{s-3\alpha-2}{s},\frac{3\alpha+2}{s})}{4^{\alpha+1}\pi\Gamma(\alpha+1)(s-3\alpha-2)}\ \left(\frac{s-3\alpha-2}{3\alpha+2}\right)^{\frac{3\alpha+2}{s}}\right)^{\frac{1}{2}}.$\\
Now suppose that we have equality in the relation (\ref{1.227}). Then we have equality in (\ref{1.14}) for $f_{r_0}$ and by means of (\ref{1.15}), we obtain
$$|f_{r_0}(x,t)|=C\ \left(1+ \ |(x,t)|^{2s} \right)^{-1},$$
and then\quad$|f(x,t)|=\left(a+b \ |(x,t)|^{2s} \right)^{-1}.$\hfill$\square$
\begin{thm}
Let $s$ be a real number such that $s>3\alpha+2$. Then for all nonzero $f\in L_\alpha^2(\mathbb{K})$ and for all measurable subset $E\subset \widehat{\mathbb{K}}$ such that $0<\gamma_\alpha(E)<+\infty$, we have
\begin{equation}\label{1.18}
 \left\|\mathcal{F}_{LB}(f)\chi_E\right\|_{\gamma_\alpha,2}< M_{\alpha,s}\ \  \gamma_\alpha(E)^\frac{1}{2}\ \left\|f\right\|^{1-\frac{3\alpha+2}{s}}_{\alpha,2}\ \left\| \ |(x,t)|^s \ f \right\|^{\frac{3\alpha+2}{s}}_{\alpha,2}
\end{equation}
    where $M_{\alpha,s}$ is the constant given by the relation $(\ref{1.12})$.
\end{thm}
\textbf{Proof}.
Suppose that the right-hand side of (\ref{1.18}) is finite. Then, according to Lemma \ref{5.2}, the function $f$ belongs to $L_\alpha^1(\mathbb{K})$ and we have
\begin{align*}
\left\|\mathcal{F}_{LB}(f)\chi_E\right\|^2_{\gamma_\alpha,2}&\leq \gamma_\alpha(E) \left\|\mathcal{F}_{LB}(f)\right\|^2_{\gamma_\alpha,\infty}\\
&\leq\gamma_\alpha(E) \left\|f\right\|^2_{\alpha,1}\\
&\leq M_{\alpha,s}^2\ \  \gamma_\alpha(E)\ \left\|f\right\|^{2-\frac{6\alpha+4}{s}}_{\alpha,2}\ \left\|\ |(x,t)|^s \ f\right\|^{\frac{6\alpha+4}{s}}_{\alpha,2}
\end{align*}
where $M_{\alpha,s}$ is the constant given by the relation (\ref{1.12}).\\
Let us prove that the equality in (\ref{1.18}) cannot hold. Suppose that there exists
a nonzero function $f\in L_\alpha^2(\mathbb{K})$
such that
$$
\left\|\mathcal{F}_{LB}(f)\chi_E\right\|^2_{\gamma_\alpha,2}=M_{\alpha,s}^2\ \  \gamma_\alpha(E)\ \left\|f\right\|^{2-\frac{6\alpha+4}{s}}_{\alpha,2}\ \left\| \ |(x,t)|^s \ f\right\|^{\frac{6\alpha+4}{s}}_{\alpha,2}
$$
Consequently, we find
\begin{equation}\label{1.19}
 \left\|\mathcal{F}_{LB}(f)\chi_E\right\|^2_{\gamma_\alpha,2}=\gamma_\alpha(E)\  \left\|\mathcal{F}_{LB}(f)\right\|^2_{\gamma_\alpha,\infty}
\end{equation}
\begin{equation}\label{1.20}
 \left\|f\right\|_{\alpha,1}=\left\|\mathcal{F}_{LB}(f)\right\|_{\gamma_\alpha,\infty},
\end{equation}
and
\begin{equation}\label{1.21}
\left\|f\right\|_{\alpha,1}=
M_{\alpha,s}\ \left\|f\right\|^{1-\frac{3\alpha+2}{s}}_{\alpha,2}\ \left\|\ |(x,t)|^s \ f\right\|^{\frac{3\alpha+2}{s}}_{\alpha,2}
\end{equation}
Applying Lemma \ref{5.2} and the relation (\ref{1.21}), we deduce that
\begin{equation}\label{1.22}
\forall (x,t)\in \mathbb{K}, \ f(x,t)= \psi(x,t)\ \left(a+b \ |(x,t)|^{2s}\right)^{-1},
\end{equation}
with $|\psi(x,t)|=1$, $a>0$ and $b>0$.\\
On the other hand, there exists $(\lambda_0,m_0)\in \hat{\mathbb{K}}$, such that
\begin{equation}\label{1.23}
\left\|\mathcal{F}_{LB}(f)\right\|_{\gamma_\alpha,\infty}=\left|\mathcal{F}_{LB}(f)(\lambda_0,m_0)\right|=e^{i\theta_0}\ \mathcal{F}_{LB}(f)(\lambda_0,m_0), \ \ \ \theta_0 \in \mathbb{R}.
 \end{equation}
 Combining now the relations (\ref{1.20}), (\ref{1.22}) and (\ref{1.23}), we get
$$\displaystyle\int_{\mathbb{K}}\left(a+b \ |(x,t)|^{2s}\right)^{-1}\left(1-e^{i\theta_0}\ \psi(x,t)\ j_{\alpha-\frac{1}{2}}(\lambda_0t)\  \mathcal{L}_{m_0}^\alpha(\lambda_0 x^2)\right)dm_\alpha(x,t)=0.$$
This implies that for almost every $(x,t)\in \mathbb{K}$,
$$e^{i\theta_0}\ \psi(x,t) \ j_{\alpha-\frac{1}{2}}(\lambda_0t)\ \mathcal{L}_{m_0}^\alpha(\lambda_0 x^2)=1.$$
Since $|\psi(x,t)|=1$, we deduce that for all $x\in \mathbb{R}_+$,
$$\left|\ j_{\alpha-\frac{1}{2}}(\lambda_0t)\ \mathcal{L}_{m_0}^\alpha(\lambda_0 x^2)\right|=1.$$
It follows that $\lambda_0=0$ and then $$\psi(x,t)=e^{-i\theta_0}.$$
Therefore the relation (\ref{1.22}) yields $$f(x,t)=e^{-i\theta_0}\ \left(a+b \ |(x,t)|^{2s}\right)^{-1};\ \ \ \ $$
Now, the relation (\ref{1.19}) means that $$\displaystyle\int_{}\displaystyle\int_{E}\left(\left\|\mathcal{F}_{LB}(f)\right\|^2_{\gamma_\alpha,\infty}-|\mathcal{F}_{LB}(f)(\lambda,m)|^2\right)d\gamma_\alpha(\lambda,m)=0$$
Hence, for almost every $(\lambda,m)\in E$, we have
\begin{equation}\label{1.24}
\left|\mathcal{F}_{LB}(f)(\lambda,m)\right|=\left\|\mathcal{F}_{LB}(f)\right\|_{\gamma_\alpha,\infty}=e^{i\theta_0}\  \mathcal{F}_{LB}(f)(0,m_0).
 \end{equation}
 $$\left|\mathcal{F}_{LB}(f)(\lambda,m)\right|=\sigma(\lambda,m)\ \mathcal{F}_{LB}(f)(\lambda,m).$$
 with $|\sigma(\lambda,m)|=1$.\\
 Then from (\ref{1.24}), for almost every $(\lambda,m)\in E$,
$$\sigma(\lambda,m)\ \mathcal{F}_{LB}(f)(\lambda,m)=e^{i\theta_0}\mathcal{F}_{LB}(f)(0,m_0),$$
and therefore $$\displaystyle\int_{\mathbb{K}}\left(a+b \ |(x,t)|^{2s}\right)^{-1}\left(1-e^{-i\theta_0}\ \sigma(\lambda,m)\ j_{\alpha-\frac{1}{2}}(\lambda t)\  \mathcal{L}_m^\alpha(\lambda x^2)\right)dm_\alpha(x,t)=0.$$
Consequently for all
$(x,t)\in \mathbb{K}$,
$$\sigma(\lambda,m)\ e^{-i\theta_0}\ j_{\alpha-\frac{1}{2}}(\lambda t)\ \mathcal{L}_m^\alpha(\lambda x^2)=1.$$
which implies that $\lambda=0$.
However, since $\gamma_\alpha(E)>0$, this contradicts the fact that for almost every $(\lambda,m)\in E$,
$$\ \mathcal{F}_{LB}(f)(\lambda,m)=\mathcal{F}_{LB}(f)(0,m_0).$$
and shows that the inequality in (\ref{1.18}) is stictly satisfied.\hfill$\square$
\begin{thm}
 Let $s=3\alpha+2$, then for all nonzero $f\in L_\alpha^2(\mathbb{K})$ and for all measurable set $E\subset \widehat{\mathbb{K}}$ such that $0<\gamma_\alpha(E)<+\infty$, we have
\begin{equation}\label{lhs}
\left\|\mathcal{F}_{LB}(f)\chi_E\right\|_{\gamma_\alpha,2}< \ C_{\alpha}\ \gamma_\alpha(E)^{\frac{1}{2(3\alpha+2)}}\  \left\|f\right\|_{\alpha,2}^{\frac{3\alpha+1}{3\alpha+2}}\left\|\ |(x,t)|^s \ f\right\|_{\alpha,2}^{\frac{1}{3\alpha+2}}
\end{equation}
with $$C_\alpha=(3\alpha+2)^2 \ (3\alpha+1)^{-\frac{1}{2(3\alpha+2)}-1}\ \  \left(\frac{B(\frac{\alpha+1}{2},\frac{2\alpha+1}{2})}{4^{\alpha+1}\pi\Gamma(\alpha+1)}\right)^{\frac{1}{2(3\alpha+2)}}$$
\end{thm}
\textbf{Proof}
$s=3\alpha+2>1 $, then using $$
\displaystyle\frac{|(x,t)|}{\varepsilon}\leq 1+\displaystyle\frac{|(x,t)|^s}{\varepsilon^{s}},\ \ \ \varepsilon>0,$$
We get $$\left\|\ |(x,t)| \ f\right\|_{\alpha,2}\leq\ \varepsilon\left\|f\right\|_{\alpha,2}+ \varepsilon^{1-s}\left\|\ |(x,t)|^s \ f\right\|_{\alpha,2}
$$
optimizing in $\varepsilon$, we obtain
$$\left\|\ |(x,t)| \ f\right\|_{\alpha,2}\leq\ s \ (s-1)^{\frac{1}{s}-1} \left\|f\right\|_{\alpha,2}^{1-\frac{1}{s}}\left\|\ |(x,t)|^s \ f\right\|_{\alpha,2}^{\frac{1}{s}}
$$
By this inequality together with $(\ref{1.1})$ taken for $s=1$, we get the result for $s=3\alpha+2$.\hfill$\square$
\bibliographystyle{unsrt}

\end{document}